\documentclass[11pt, final]{amsart}
\usepackage{amssymb}

\theoremstyle{plain}
\newtheorem{theorem}{Theorem}
\newtheorem{proposition}[theorem]{Proposition}

\newtheorem{corollary}[theorem]{Corollary}

\theoremstyle{remark}

\def\A{\mathcal A}

\def\F{\mathbb F}
\def\K{\mathbb K}

\def\MnF{M_n(\mathbb F)}
\def\N{\mathbb N}
\def\ot{\otimes}
\def\R{\mathbb R}

\def\rank{\operatorname{rank}}
\def\si{\operatorname{eigv}}
\def\tr{\operatorname{tr}}

\begin{document}
\title[]{Local automorphisms of operator algebras on Banach spaces}
\author{LAJOS MOLN\'AR}
\address{Institute of Mathematics and Informatics\\
         University of Debrecen\\
         4010 Debrecen, P.O.Box 12, Hungary}
\email{molnarl@math.klte.hu}
\thanks{  This research was supported by the
          Hungarian National Foundation for Scientific Research
          (OTKA), Grant No. T030082, T031995, and by
          the Ministry of Education, Hungary, Reg.
          No. FKFP 0349/2000}
\date{\today}
%\date{{\LARGE \today}}
\begin{abstract}
In this paper we extend a result of \v Semrl stating that
every 2-local automorphism of the full operator algebra on a separable
infinite dimensional Hilbert space is an automorphism. In fact,
besides separable Hilbert spaces, we
obtain the same conclusion for the much larger class of Banach spaces with
Schauder bases. The proof rests on an analogous statement concerning
the 2-local automorphisms of matrix algebras of which statement we
present a short proof. The need
to get such a proof was formulated in \v Semrl's paper.
\end{abstract}
\maketitle

\section{Introduction and Statement of the Results}

The importance of the study of automorphism groups of algebraic
structures needs no justification. Recently, we have published
several results on the local behaviour of the automorphism groups
of operator algebras. Roughly speaking, we have studied the
question of how the automorphisms of certain operator algebras (or
of function algebras) are determined by their local actions. One
area of such investigations was initiated by P. \v Semrl who
introduced the concept of 2-local automorphisms. If $\A$ is any
algebra, then a transformation (no linearity is assumed) $\phi:\A
\to \A$ is called a 2-local automorphism of $\A$ if for every
$A,B\in \A$ there exists an (algebra) automorphism $\phi_{A,B}$ of
$\A$ such that $\phi(A)=\phi_{A,B}(A)$ and
$\phi(B)=\phi_{A,B}(B)$. It is then a remarkable fact on the
algebra $\A$ if every 2-local automorphism of $\A$ is an
automorphism. Observe that for 1-local automorphisms (the
definition should be self-explanatory) without any extra
properties reasonable results of that kind can not be expected.

The first result on 2-local automorphisms is due to \v Semrl
\cite{Semrl2} who proved that if $H$ is a separable infinite
dimensional Hilbert space, then every 2-local automorphism of the
algebra $B(H)$ of all bounded linear operators on $H$ is an
automorphism. Moreover, since in the finite dimensional case he
was only able to get a long proof involving tedious computations,
\v Semrl suggested it would be of some interest to produce a short
proof in that case. One aim of this paper is to present such a
proof. Furthermore, we extend \v Semrl's result quite
significantly, namely, in addition to the case of separable
Hilbert spaces, we get the same conclusion for general Banach
spaces with Schauder bases. We recall that most classical
separable Banach spaces have Schauder bases (see the first chapter
in \cite{Lindenstrauss}) and hence our generalization really has
sense.

Now we turn to the formulation of our assertions. If $\F$ is a
field, then $M_n(\F)$ stands for the algebra of all $n\times n$
matrices with entries in $\F$. It is well-known that the
automorphisms of the algebra $\MnF$ are exactly the
transformations $A\mapsto TAT^{-1}$ where $T\in \MnF$ is
nonsingular (see, for example, \cite[Lemma, p. 230]{Pierce}). In
what follows, if $A\in \MnF$, we denote by $\si(A)$ the system of
all eigenvalues of $A$ listed according to multiplicity. We
emphasize that $\si(A)$ generally differs from the spectrum of $A$
as a linear operator.

Our first theorem gives a nonlinear characterization of the
automorphisms of $\MnF$.

\begin{theorem}\label{T:locbanmat}
Suppose that $\F$ is an algebraically closed field of characteristic 0.
Let $\phi:\MnF \to \MnF$ be a transformation (linearity is not assumed)
such that
\begin{equation}\label{E:locban1}
\si(\phi(A)\phi(B))=\si(AB)
\end{equation}
holds for every $A,B\in \MnF$. Then there exists a nonsingular matrix
$T\in \MnF$ and $\lambda \in \{-1, 1\}$ such that
$\phi$ is either of the form
\begin{equation*}
\phi(A)=\lambda TAT^{-1} \qquad (A\in \MnF)
\end{equation*}
or of the form
\[
\phi(A)=\lambda TA^tT^{-1} \qquad (A\in \MnF).
\]
\end{theorem}

We mention that similar characterizations of the automorphisms of
operator algebras and of function algebras can be found in \cite{Molnar}
but there we had to assume that the transformations in
question are all surjective. The main advantage of the present theorem
is that we can omit that assumption and this fact makes us possible to
get the following
corollary concerning 2-local automorphisms of matrix algebras over
general fields.

\begin{corollary}\label{C:locbanmat}
Let $\F$ be an algebraically closed field of characteristic 0. Then
every 2-local automorphism of $\MnF$ is an automorphism.
\end{corollary}

Since the field $\R$ of real numbers is not algebraically closed, we
formulate the analogous statement for $M_n(\R)$
as a separate proposition.

\begin{proposition}\label{P:locbanmat}
Every 2-local automorphism of $M_n(\R)$ is an automorphism.
\end{proposition}

The arguments applied in the proof of our results above lead us
to an extension of \v Semrl's theorem for the case of certain
subalgebras of the algebra $B(X)$ of all
bounded linear operators on the Banach space $X$ with a Schauder basis.
We believe that with some extra work \v Semrl's
original proof could also be improved to produce the following result
(although only in the infinite dimensional case) but our
proof has the advantage that it puts both the finite and the infinite
dimensional cases into the same perspective.

\begin{theorem}\label{T:locbanop}
Let $X$ be a real or complex Banach space with a Schauder basis.
Suppose that $\A$ is a subalgebra of $B(X)$ which contains the ideal of
all compact
operators on $X$. If $\phi$ is a 2-local automorphism of $\A$, then
$\phi$ is an automorphism of $\A$.
\end{theorem}

\section{Proofs}

In the proof of our first result we shall need the following
folklore result whose proof requires only elementary linear
algebra. If $X$ is a linear space over $\F$ and $A,B: X\to X$ are
linear operators of rank at least 2 with the property that for
every $x\in X$ the vectors $Ax,Bx$ are linearly dependent, then
the operators $A,B$ are linearly dependent.

\begin{proof}[Proof of Theorem~\ref{T:locbanmat}]
Since the statement is obvious for $n=1$, we suppose that $n\geq
2$. Denote by $\tr$ the usual trace functional on $\MnF$. It
follows from the property \eqref{E:locban1}, i.e., from the
equality $\si(\phi(A)\phi(B))=\si(AB)$ that
\begin{equation}\label{E:locban2}
\tr \phi(A)\phi(B)=\tr AB
\end{equation}
holds for every $A,B\in \MnF$.

As usual, denote by $E_{ij}\in \MnF$
the matrix whose $ij$ entry is 1 and its all other entries are 0.
We assert that the $\phi(E_{ij})$'s are linearly independent.
Suppose that
\[
\sum_{i,j} \lambda_{ij}\phi(E_{ij}) =0
\]
for some scalars $\lambda_{ij}\in \mathbb F$.
Fix indices $k,l\in \{ 1, \ldots, n\}$. We have
\[
\sum_{i,j} \lambda_{ij}\phi(E_{ij})\phi(E_{kl}) =0.
\]
Taking trace, we obtain
\[
\sum_{i,j} \lambda_{ij}\tr \phi(E_{ij})\phi(E_{kl}) =0.
\]
By the property \eqref{E:locban2} of $\phi$, it follows that
\[
\sum_{i,j} \lambda_{ij}\tr E_{ij}E_{kl} =0.
\]
Since $E_{ij}E_{kl}=\delta_{jk} E_{il}$, from this equality we easily
deduce that $\lambda_{lk}=0$. As $k,l$ were arbitrary, it follows that
the matrices $\phi(E_{ij})$, $i,j\in \{ 1, \ldots, n\}$ form a linearly
independent set in $\MnF$.
This implies that the range of $\phi$ linearly generates $\MnF$.

We are now in a position to prove that $\phi$ is linear.
If $A,B \in \MnF$, we compute
\begin{equation*}
\begin{gathered}
\tr \phi(A+B)\phi(C)=\tr (A+B)C=
\tr AC +\tr BC=\\
\tr (\phi(A)\phi(C)+\phi(B)\phi(C))=
\tr (\phi(A)+\phi(B))\phi(C)
\end{gathered}
\end{equation*}
for any $C\in \MnF$. Since the linear span of the $\phi(C)$'s is the
whole space $\MnF$, it follows that $\phi(A+B)=\phi(A)+\phi(B)$.
The homogenity of $\phi$ can be proved in a very similar way.
So, $\phi$ is a surjective linear transformation on $\MnF$ which is
hence bijective.

Let $A\in \MnF$ be of rank one. Then, using \eqref{E:locban1} and the
surjectivity of $\phi$, it follows
that  $\phi(A)B, B\phi(A)$ both have at most one nonzero eigenvalue and
its multiplicity is one. This implies that $\phi(A)$ is of rank one.
Consequently, $\phi$ preserves the rank-one elements of $\MnF$.
The form of all such linear transformations on $\MnF$ is well-known. It
follows from
\cite{Marcus} that there are nonsingular matrices $T,S$ in $\MnF$ such
that $\phi$ is either of the form
\begin{equation}\label{E:locban3}
\phi(A)=TAS \qquad (A\in \MnF)
\end{equation}
or of the form
\[
\phi(A)=TA^tS \qquad (A\in \MnF).
\]
Without loss of generality we can assume that $\phi$ is of the first
form.
Composing our transformation $\phi$ with the automorphism $A\mapsto
SAS^{-1}$, it follows that we can assume that the matrix $S$ appearing
in \eqref{E:locban3} is the identity $I$.
We have
\begin{equation}\label{E:locban4}
\tr TATB=\tr AB
\end{equation}
for every $A,B\in\MnF$. If we fix $A\in \MnF$ for a moment and let $B$
run through $\MnF$, it follows from \eqref{E:locban4} that
$TAT=A$.
Now, if $A$ runs through the set of all rank-one matrices in $\MnF$,
we deduce from the equality $TAT=A$
that $T$, as a linear operator on the linear space $X$, has the
property that
for every $x\in X$, the vector $Tx$ is a scalar multiple of $x$. This
gives
us that $T$ is of the form $\lambda I$ for some scalar $\lambda\in \F$.
Clearly, we have $\lambda^2=1$ and this completes the proof of our
statement.
\end{proof}

\begin{proof}[Proof of Corollary~\ref{C:locbanmat}]
First we recall again that the automorphisms of $\MnF$ are exactly
the transformations of the form $A\mapsto TAT^{-1}$ with some
nonsingular matrix $T\in \MnF$.
Let now $\phi:\MnF \to \MnF$ be a 2-local automorphism.
Clearly, $\phi$ satisfies \eqref{E:locban1}.
Since $\phi(I)=I$, we infer from Theorem~\ref{T:locbanmat} that there
exists a nonsingular matrix $T\in \MnF$ such that $\phi$ is either of
the form
\[
\phi(A)=TAT^{-1} \qquad (A\in \MnF)
\]
or of the form
\[
\phi(A)=TA^tT^{-1} \qquad (A\in \MnF).
\]
Choosing $A,B\in \MnF$ such that $AB=0, BA\neq 0$, it
follows from the
2-local property of $\phi$ that this second possibility above
is excluded and we obtain the desired conclusion.
\end{proof}

\begin{proof}[Proof of Proposition~\ref{P:locbanmat}]
The statement follows from obvious modifications of the proof of
Corollary~\ref{C:locbanmat} and that of Theorem~\ref{T:locbanmat}.
The only difference is that instead of referring to the
result of Marcus and Moyls on the form of all linear
transformations on $\MnF$ preserving rank-one matrices, here we have to
refer to
an analogous statement concerning the algebra of matrices over the real
field due to Omladi\v c and \v Semrl \cite[Main Theorem]{OmladicSemrl}.
\end{proof}

In the proof of Theorem~\ref{T:locbanop}
we need the following notation and
definitions.
Let $X$ be a (real or complex) Banach space. The algebra of all bounded
linear operators on $X$ is denoted by $B(X)$ and $F(X)$ stands for
the ideal of all finite rank operators
in $B(X)$.
An operator $P\in B(X)$ is called an idempotent if $P^2=P$.
Two idempotents $P,Q$ in $B(X)$ are said to be
(algebraically) orthogonal if $PQ=QP=0$.
The dual space of $X$ (that is the set of all bounded linear functionals
on $X$) is denoted by $X'$. The Banach space adjoint of
an operator $A\in B(X)$ is denoted by $A'$.
If $x\in X$ and $f\in X'$, then $x\ot f$ stands for the
operator (of rank at most 1) defined by
\[
(x\ot f)(z)=f(z)x \qquad (z\in X).
\]
It requires only elementary computation to show that
\[
A \cdot x\ot f=(Ax) \ot f,
\quad
x\ot f \cdot A=x\ot (A'f),
\quad
x\ot f \cdot y\ot g= f(y) x\ot g
\]
hold for every $x,y\in X$, $f,g\in X'$ and $A\in B(X)$.
It is easy to see that the elements of $F(X)$ are exactly the operators
$A\in B(X)$ which can be written as a finite sum
\begin{equation}\label{E:zen5}
A=\sum_i x_i \ot f_i
\end{equation}
for some $x_1, \ldots, x_n\in X$ and $f_1, \ldots, f_n\in X'$.
Using this representation, the trace of $A$ is defined by
\[
\tr A=\sum_i  f_i(x_i).
\]
It is known that $\tr A$ is well-defined, that is, it does not depend on
the particular representation \eqref{E:zen5} of $A$. It is easy to see
that $\tr$ is a linear functional on $F(X)$ with the property that
\[
\tr TA=\tr AT
\]
holds for every $A\in F(X)$ and $T\in B(X)$.
Finally, we recall that, similarly to the
case of matrix algebras, every automorphism of
any subalgebra of $B(X)$ which contains $F(X)$ (these are the so-called
standard operator algebras on $X$) is of the form $A \mapsto TAT^{-1}$
where $T\in B(X)$ is invertible (see, for example, \cite{Semrl1}).

\begin{proof}[Proof of Theorem~\ref{T:locbanop}]
In view of our previous results we can assume that $X$ is infinite
dimensional.
Let $\phi$ be a 2-local automorphism of $\A$.
If $P\in \A$ is a
finite rank idempotent, then $\phi(P)=\tilde P$ is also a finite rank
idempotent and the ranks of $P$ and $\tilde P$ are the same.
Consider the subalgebra $\A_P$ of $\A$ which consists of all
operators $A\in \A$ for which $PAP=A$. It follows from the 2-local
property
of $\phi$ that $\phi(P)\phi(A)\phi(P)=\phi(A)$. Consequently, $\phi$
maps $\A_P$ into $\A_{\tilde P}$. Clearly, both algebras
$\A_P$ and $\A_{\tilde P}$ are isomorphic to $M_n(\K)$ ($\K$ stands for
the real or complex field) where $n=\rank P=\rank \tilde P$.
It is easy to see that $\phi$ has the property
\eqref{E:locban2}. Now, one can follow the arguments given in
the proofs of our previous results to show that $\phi$ is
linear and multiplicative on $\A_P$.
Since $P\in \A$ was an arbitrary finite rank idempotent, we can infer
that
the restriction $\psi$ of $\phi$ onto $F(X)$ is an algebra endomorphism
of $F(X)$ which preserves the rank.

Since $\psi$ is an algebra homomorphism of $F(X)$, we have an injective
linear operator $T:X \to X$ such that
\begin{equation}\label{E:zen1}
TA=\psi(A)T \qquad (A\in F(X)).
\end{equation}
Indeed, similarly as in \cite{Semrl1} we define
\[
Tx=\psi(x \ot f_0)y_0 \qquad (x\in X)
\]
where $y_0\in X$ and $f_0\in X'$ are fixed such that $\psi(x\ot
f_0)y_0\neq 0$ for some $x\in X$. It is evident that $T$ is a nonzero
linear operator on $X$. It follows from the multiplicativity of $\psi$
that $TA=\psi(A)T$ $(A\in F(X))$. To see
the injectivity of $T$, let $Tx=0$ and $x\neq 0$. Then we have
$TAx=\psi(A)Tx=0$ for every $A\in F(X)$. This gives us that $Ty=0$
holds for every $y\in X$ which is an obvious contradiction.

We show that $T$ is continuous. To see this, we apply the closed graph
theorem. Let $(x_n)$ be a sequence in $X$ and $y\in X$ be such that
$x_n \to 0$ and $Tx_n \to y$. We have to show that $y=0$.
Since $TA$ is continuous for every $A\in F(X)$, from \eqref{E:zen1} we
infer that
\begin{equation}\label{E:br5}
\psi(A)y=0 \qquad (A\in F(X)).
\end{equation}

Here is the point where we have to use that $X$ has a Schauder basis.
Namely, this condition implies that
we have a sequence $(P_n)$ of pairwise orthogonal rank-one
idempotents in $B(X)$ whose sum strongly converges to $I$.
In particular, $(P_n)$ is uniformly bounded.
Clearly, we can write $P_n=x_n \ot f_n$ where $x_n \in X, f_n \in X'$
are such that $f_i(x_j)=\delta_{ij}$ (the Kronecker symbol) and $\|
x_i\|
<m, \| f_j\| <M$ $(i,j\in \N)$ for some positive real numbers $m,M$.
Choose a sequence $(\lambda_n)$ of
positive real numbers with the property that
$\sum_{n=k+1}^\infty
\lambda_n <\lambda_k$ holds for every $k\in \N$. For example, one can
define $\lambda_n =(1/3)^n$ $(n\in \N)$. In particular, it follows that
$\sum_n \lambda_n P_n$ converges in the norm topology and hence
its sum is a compact operator.
Pick another sequence $(\mu_n)$ of positive real numbers for which
$\sum_n
\mu_n <\infty$. The operator $\sum_n \mu_n x_n \ot f_{n+1}$ is also
compact.

By the 2-local property of $\phi$, composing $\phi$ with an automorphism
of $\A$ if necessary, we can (and do) assume that for
the particular operators $\sum_n \lambda_n P_n$
and $\sum_n \mu_n x_n\ot f_{n+1}$
we have
\begin{equation}\label{E:br2}
\phi(\sum_n \lambda_n P_n)=\sum_n \lambda_n P_n
\text{ and }
\phi(\sum_n \mu_n x_n \ot f_{n+1})=\sum_n \mu_n x_n \ot f_{n+1}.
\end{equation}

Let $n_0\in \N$ be arbitrary. By the 2-local property of $\phi$ we
have an invertible bounded linear operator $U\in B(X)$ such that
\[
\phi(\sum_n \lambda_n P_n)=U(\sum_n \lambda_n P_n)U^{-1}
\quad \text{and} \quad
\phi(P_{n_0})=UP_{n_0}U^{-1}.
\]
From
\[
\sum_n \lambda_n UP_nU^{-1}=
U(\sum_n \lambda_n P_n)U^{-1}=
\phi(\sum_n \lambda_n P_n)=
\sum_n \lambda_n P_n
\]
we can infer that $UP_nU^{-1}=P_n$ holds for every $n\in \N$.
Indeed, let $UP_nU^{-1}=Q_n$ $(n\in \N)$.
Dividing both sides of the equality
\[
\sum_n \lambda_n Q_n =
\sum_n \lambda_n P_n
\]
by $\lambda_1$ and then taking the $k$th powers of the operators on
both sides and, finally,
letting $k$ tend to infinity, by the property of the sequence
$(\lambda_n)$ we easily obtain that $Q_1=P_1$. Therefore, we have
\[
\sum_{n=2}^\infty \lambda_n Q_n =
\sum_{n=2}^\infty \lambda_n P_n
\]
and one can proceed in the same way to show that $Q_n=P_n$ holds for
every
$n\in \N$. In particular, we have $\phi(P_{n_0})=Q_{n_0}=P_{n_0}$. But
$n_0$
was arbitrary and hence we obtain that
\begin{equation}\label{E:br1}
\phi(P_n)=P_n  \qquad (n\in \N).
\end{equation}

We now go back to the proof that $T$ is continuous.
By \eqref{E:br5} we have $\psi(P_n)y=0$ $(n\in \N)$ and this,
together with \eqref{E:br1}, yields that $P_ny
=0$ holds for every $n\in \N$. This implies that $y=0$ verifying the
continuity of $T$.

Beside the fact that $\psi$ is an algebra endomorphism of $F(X)$, we
know that $\psi$ is rank preserving. The form of
linear rank preservers on operator algebras is known.
In fact, a description of such transformations
is given, for example, in \cite{Hou}. It follows from the arguments
used in the first
half of Section I in \cite{Hou} that we have two possibilities:
either there are linear operators $S:X\to X$ and $R:X' \to X'$ such that
\[
\psi(x\ot f)=Sx\ot Rf \qquad (x\in X, f\in X')
\]
or there are linear operators $R:X\to X'$ and $S:X' \to X$ such that
\[
\psi(x\ot f)=Sf\ot Rx \qquad (x\in X, f\in X').
\]
(To be honest we note that in \cite{Hou} only complex spaces were
considered but the arguments given there also apply to the real case.)
This second possibility can be excluded easily after referring to
the multiplicativity of $\psi$.
So, we have linear operators $S:X\to X$ and $R:X' \to X'$ such that
\begin{equation}\label{E:zen3}
\psi(x\ot f)=Sx\ot Rf \qquad (x\in X, f\in X').
\end{equation}
It follows from \eqref{E:zen1} that
\[
Tx\ot f =
T \cdot x\ot f =Sx\ot Rf \cdot T=
Sx\ot T'Rf.
\]
This implies that $Tx,Sx$ are linearly dependent for every $x\in X$.
Hence, $S$ is a scalar multiple of $T$.
Therefore, in what follows we can (and do) assume that the linear
operator $S$ appearing in \eqref{E:zen3} is equal to $T$.

By the local form of $\phi$ it follows that $\psi$ is trace preserving.
Consequently, we have
\[
(Rf)(Tx)=f(x) \qquad (x\in X, f\in X').
\]
This means that $T'R$ equals the identity on $X'$. In particular,
the range of $T'$ is closed which is well-known to imply that the
range of $T$ is also closed. On the other hand, $T$ has dense
range which follows from \eqref{E:zen1} and \eqref{E:br1}.
Therefore, we can infer that $T,T'$ are invertible and
$R={T'}^{-1}={T^{-1}}'$. This gives us that
\begin{equation}\label{E:zen30}
\psi(x\ot f)=T \cdot x\ot f \cdot T^{-1} \qquad (x\in X, f\in X')
\end{equation}
and hence we have $\phi(A)=\psi(A)=TAT^{-1}$ for every $A\in F(X)$.
We show that $\phi(A)=TAT^{-1}$ holds also for every $A\in \A$.
In order to do this,
let $A\in \A$ be arbitrary. Pick any $x\in X, f\in X'$ and let $B=x\ot
f$. By the 2-local property of $\phi$ we have an invertible element $U$
of $B(X)$ for which
\[
\phi(A)=UAU^{-1}
\quad \text{and} \quad
\phi(B)=UBU^{-1}.
\]
Since $B\in F(X)$, we can compute
\begin{equation*}
\begin{gathered}
TBT^{-1}\phi(A)TBT^{-1}=\phi(B)\phi(A)\phi(B)=\\
UBU^{-1} UAU^{-1} UBU^{-1}= UBABU^{-1}.
\end{gathered}
\end{equation*}
Taking traces in this equality, we get
\[
f(T^{-1}\phi(A)Tx)f(x)=f(Ax)f(x).
\]
This implies that for any $x\in X, f\in X'$ with $f(x)\neq 0$ we have
\begin{equation}\label{E:zen7}
f(T^{-1}\phi(A)Tx)=f(Ax).
\end{equation}
Moreover,
if $x\in X$, $0\neq f\in X'$ and $f(x)=0$, then we can find a sequence
$(x_n)$ in $X$ such
that $x_n \to x$ and $f(x_n)\neq 0$ for every $n\in \N$. Hence, we
deduce that the equality \eqref{E:zen7}
holds without any assumption on $x$ and $f$. This gives us that
\begin{equation}\label{E:br3}
\phi(A)=TAT^{-1} \qquad (A\in \A).
\end{equation}
To complete the proof it remains to show that in fact we have
$\phi(A)=A$ for every $A\in \A$ (this is to prove that $\phi$ is
surjective).
From \eqref{E:br1} we learn that $\phi(x_n \ot f_n)=x_n \ot f_n$ for
every $n\in \N$.
Due to \eqref{E:zen30} this readily imples that $Tx_n=\epsilon_n x_n$
and ${T^{-1}}'f_n =\delta_n f_n$ hold for some scalars $\epsilon_n,
\delta_n$. Moreover, we have $\epsilon_n \delta_n=1$ for every $n\in
\N$. Obviously, we can suppose that $\epsilon_1=\delta_1=1$.
We have supposed in \eqref{E:br2} that
\[
\phi(\sum_n \mu_n x_n \ot f_{n+1})=
\sum_n \mu_n x_n \ot f_{n+1}.
\]
By \eqref{E:br3} this yields
\[
\sum_n \mu_n Tx_n \ot {T^{-1}}'f_{n+1}=
\sum_n \mu_n x_n \ot f_{n+1}.
\]
Therefore, we have
\[
\sum_n \epsilon_n\delta_{n+1} \mu_n x_n \ot f_{n+1}=
\sum_n \mu_n x_n \ot f_{n+1}.
\]
Considering the values of the operators on both sides at $x_2, x_3,
\ldots$ one after the other,
we get in turn that $\delta_2=1, \epsilon_2=1,\enskip \delta_3=1,
\epsilon_3=1, \enskip
\ldots$. Therefore, $T$ is the identity on $X$ and we have
\[
\phi(A)=A \qquad (A\in \A).
\]
This completes the proof.
\end{proof}

% Bibliography
\bibliographystyle{amsplain}

\end{document}